\documentclass[12pt, leqno]{amsart}
\usepackage{lscape}
\usepackage{graphicx}
\usepackage{amssymb,amsmath,amsthm,amscd}
\usepackage{mathrsfs}
\usepackage{enumerate}
\usepackage{verbatim}
\usepackage[usenames,dvipsnames]{color}
\usepackage[colorlinks=true, pdfstartview=FitV,  linkcolor=blue,citecolor=blue,urlcolor=blue, bookmarks=false]{hyperref}

\usepackage[all]{xy}
\allowdisplaybreaks[4]
\usepackage{oldgerm}

\newcommand{\nc}{\newcommand}
\nc{\cmt}{\fbox{comment}}

\newcommand{\bj}{\begin{jaune}}
\newcommand{\ej}{\end{jaune}}

\theoremstyle{plain}
\newtheorem*{theorem*}{Theorem}
\newtheorem{lemma}{Lemma}%[section]
\newtheorem{prop}[lemma]{Proposition}
\newtheorem{theorem}[lemma]{Theorem}

\newcommand{\Prop}{\begin{prop}}
\newcommand{\enprop}{\end{prop}}
\newcommand{\Lemma}{\begin{lemma}}
\newcommand{\enlemma}{\end{lemma}}
\newcommand{\Th}{\begin{theorem}}
\newcommand{\enth}{\end{theorem}}
\newtheorem{corollary}[lemma]{Corollary}
\newcommand{\Cor}{\begin{corollary}}
\newcommand{\encor}{\end{corollary}}
\newtheorem{definition}[lemma]{Definition}
\newtheorem*{conjecture}{Conjecture}
\newcommand{\Def}{\begin{definition}}
\newcommand{\edf}{\end{definition}}
\newtheorem{sublemma}[lemma]{Sublemma}
\newcommand{\Sub}{\begin{sublemma}}
\newcommand{\ensub}{\end{sublemma}}

\theoremstyle{definition}
\newtheorem*{remark}{Remark}

\newtheorem{Convention}[lemma]{Convention}
\newcommand{\Conv}{\begin{Convention}}
\newcommand{\enconv}{\end{Convention}}
\nc{\Con}{\begin{conjecture}}
\nc{\encon}{\end{conjecture}}
\nc{\Rem}{\begin{remark}}
\nc{\enrem}{\end{remark}}

\newcommand{\C}{{\mathbb C}}
\newcommand{\Q}{\mathbb {Q}}

\newcommand{\Z}{{\mathbb Z}}
\newcommand{\B}{{\mathbf{B}}}
\newcommand{\A}{{\mathbb Z [q^{\pm1}]}}

\newcommand{\gl}{{\mathfrak{gl}}}

\newcommand{\seteq}{\mathbin{:=}}

\newcommand{\g}{\mathfrak{g}}

\newcommand{\n}{\mathfrak{n}}

\newcommand{\Uq}[1][{\mathfrak{g}}]{{U_q(#1)}}
\newcommand{\Uqm}[1][{\mathfrak{g}}]{{U_q^-(#1)}}

\newcommand{\Hom}{\operatorname{Hom}}
\newcommand{\End}{\operatorname{End}}

\newcommand{\isoto}[1][]{\mathop{\xrightarrow%
[{\raisebox{.3ex}[0ex][.3ex]{$\scriptstyle{#1}$}}]%
{{\raisebox{-.6ex}[0ex][-.6ex]{$\mspace{2mu}\sim\mspace{2mu}$}}}}}
\nc{\N}{\Z_{\ge0}}
\newcommand{\eq}{\begin{eqnarray}}
\newcommand{\eneq}{\end{eqnarray}}

\newcommand{\hs}{\hspace*}

\newcommand{\To}[1][{\hs{2ex}}]{\xrightarrow{\,#1\,}}

\newcommand{\eqn}{\begin{eqnarray*}}
\newcommand{\eneqn}{\end{eqnarray*}}
\newcommand{\on}{\operatorname}

\newcommand{\bna}{\be[{\rm(a)}]}

\newcommand{\QED}{\end{proof}}
\newcommand{\Proof}{\begin{proof}}

\newcommand{\soplus}{\mathop{\mbox{\normalsize$\bigoplus$}}\limits}

\newcommand{\id}{\on{id}}
\newcommand{\ba}{\begin{array}}
\newcommand{\ea}{\end{array}}

\newcommand{\set}[2]{\left\{#1 \mid #2 \right\}}

\newcommand{\eqsub}{\begin{subequations}\begin{eqnarray}}
\newcommand{\eneqsub}{\end{eqnarray}\end{subequations}}

\newcommand{\ol}{\overline}

\nc{\la}{\lambda}
\nc{\lam}{\lambda}
\nc{\U}[1][\g]{U_q(#1)}
\nc{\te}{\tilde{e}}
\nc{\tei}{\tilde{e}_i}
\nc{\tf}{\tilde{f}}
\nc{\tfi}{\tilde{f}_i}
\nc{\tU}{\widetilde U_q(\g)}
\nc{\tE}{\tilde{E}}
\nc{\tF}{\widetilde{F}}
\nc{\tK}{\widetilde{K}}

\nc{\tk}{\tilde{k}}
\nc{\tkone}{\tk_{\ol{1}}}
\nc{\teone}{\tilde{e}_{\ol{1}}}
\nc{\tfone}{\tilde{f}_{\ol{1}}}

\nc{\teibar}{\tilde{e}_{\ol{i}}} \nc{\tfibar}{\tilde{f}_{\ol{i}}}
\nc{\tki}{{\tk}_{\ol {i}}}

\nc{\BZ}{{\mathbb{Z}}}
\nc{\al}{\alpha}
\nc{\qs}{{q}}
\nc{\lan}{\langle}
\nc{\ran}{\rangle}
\nc{\re}{{\mathrm{re}}}
\nc{\wt}{\operatorname{wt}}
\nc{\ch}{\operatorname{ch}}
\nc{\Um}[1][\g]{U^-_q(#1)}
\nc{\Ue}{U^+_q(\g)}
\nc{\eps}{\varepsilon}
\nc{\vphi}{\varphi}
\nc{\sphi}{\varphi^*}
\nc{\seps}{\varepsilon^*}
\newcommand{\bB}{\mathbf{B}}
\nc{\nn}{\nonumber}
\def\max{{\mathop{\mathrm{max}}}}
\nc{\vph}{\varphi}
\nc{\cls}{{\operatorname{cl}}}
\nc{\Wt}{{\operatorname{Wt}}}
\nc{\Us}{U'_q(\g)}
\nc{\La}{\Lambda}
\nc{\tLa}{\widetilde\Lambda}
\nc{\ro}{{\rm(}}
\nc{\rf}{{\rm)}}
\nc{\norm}{{\mathrm{norm}}}
\nc{\qbox}{\quad\mbox}
\nc{\braid}{{\mathfrak{B}}}
\nc{\Ad}{\operatorname{Ad}}
\nc{\Aut}{\operatorname{Aut}}
\nc{\dt}[1]{\tilde{\tilde #1}}
\nc{\Sn}{S^{{\mathrm{norm}}}}
\nc{\aff}{{\rm{aff}}}
\nc{\rk}{{\mathrm{rk}}}
\nc{\tP}{\widetilde{P}}
\nc{\tW}{\widetilde{W}}
\nc{\Dyn}{\mathrm{Dyn}}
\nc{\tD}{\widetilde{\Delta}}
\nc{\height}[1]{{\operatorname{ht}}(#1)}
\nc{\bl}{\bigl(}
\nc{\br}{\bigr)}
\nc{\Hecke}{\mathrm{H}}
\nc{\HA}{\Hecke^{\mathrm{A}}}
\nc{\HB}{\Hecke^{\mathrm{B}}}
\newcommand{\scbul}{{\,\raise1pt\hbox{$\scriptscriptstyle\bullet$}\,}}
\nc{\vac}{{\phi}}
\nc{\Bt}{\B_\theta(\g)}
\nc{\be}{\begin{enumerate}}
\nc{\ee}{\end{enumerate}}
\nc{\low}{{\mathrm{low}}}
\nc{\upper}{{\mathrm{up}}}
\nc{\Zodd}{\Z_{\mathrm{odd}}}
\nc{\Ft}[1][n]{\mathbb{P}\mathrm{ol}_{#1}}
\nc{\Ftf}[1][n]{\widetilde{\mathbb{P}\mathrm{ol}}_{#1}}
\nc{\KA}{\on{K}^{\mathrm{A}}}
\nc{\KB}{\on{K}^{\mathrm{B}}}
\nc{\Res}{\on{Res}}
\nc{\Fc}[1][{n,m}]{\mathbf{F}_{#1}}
\nc{\tphi}{\tilde{\varphi}}
\nc{\CO}{\mathscr{O}}
\nc{\inte}{\mathrm{int}}

\nc{\vs}{\vspace*}
\nc{\tLt}{\widetilde{L}}
\nc{\tL}{\widetilde{\Lambda}}
\nc{\tu}{\tilde{u}}
\nc{\noi}{\noindent}
\nc{\heigh}{\mathfrak{t}}
\nc{\lowest}{\mathfrak{l}}
\nc{\rootl}{\mathsf{Q}}
\nc{\cl}{\colon}
\nc{\uqpg}{U'_q(\mathfrak g)}
\nc{\uq}{\uqpg}
\nc{\Oh}{\widehat{\mathcal{O}}}

\nc{\KLR}{KLR algebra}
\nc{\KLRs}{KLR algebras}
\nc{\cor}{\mathbf{k}}
\nc{\cora}{{\cor(A)}}
\nc{\haut}{\mathrm{ht}}
\nc{\tens}{\mathop\otimes}
\nc{\gmod}{\mbox{-$\mathrm{gmod}$}}
\nc{\gMod}{\mbox{-$\mathrm{gMod}$}}
\nc{\proj}{\mbox{-$\mathrm{proj}$}}
\nc{\gproj}{\mbox{-$\mathrm{gproj}$}}
\nc{\smod}{\mbox{-$\mathrm{mod}$}}
\nc{\Mod}{\mbox{-$\mathrm{Mod}$}}
\nc{\Rnorm}{R^{\rm{norm}}}

\nc{\Vhat}{\widehat{V}}
\nc{\F}{\mathcal{F}}

\def\T{{\mathcal T}}

\nc{\fd}[1][A]{\on{\mathrm{flat.dim}_{#1}}}
\nc{\bP}{{\mathbb{P}}}
\nc{\bPh}{\widehat{\mathbb{P}}}
\nc{\bK}[1][{n}]{\widehat{\mathbb{K}}_{#1}}
\nc{\bV}[1][{n}]{\widehat{V}^{\otimes{#1}}}
\nc{\bVK}[1][{n}]{\widehat{V}^{\otimes{#1}}_{\widehat{\mathbb{K}}}}
\nc{\hV}{\widehat{V}}
\nc{\opp}{\mathrm{opp}}
\nc{\col}{\colon}
\nc{\bnum}{\be[{\rm(i)}]}
\nc{\bnam}{\be[{\rm(a)}]}
\nc{\oep}{\epsilon}
\nc{\qtext}[1][{and}]{\quad\text{#1}\quad}

\nc{\qtextq}[1]{\quad\text{#1}\quad}
\nc{\longtwoheadrightarrow}[1][]{\xymatrix{\ar@{->>}[r]^-{{#1}}&}}
\nc{\epiTo}[1][]{\longtwoheadrightarrow[{#1}]}
\nc{\epito}{\twoheadrightarrow}
\nc{\monoTo}[1][]{\xymatrix{\ar@{>->}[r]^-{{#1}}&}}
\nc{\sym}{\mathfrak{S}}
\nc{\inp}[1]{{({#1})_{\mathrm{n}}}}
\nc{\rtl}{\rootl}
\nc{\wtd}{\widetilde}
\nc{\etens}{\boxtimes}
\nc{\ds}[1]{\mathrm{d}(#1)}
\nc{\rmat}[1]{{\mathbf{r}}_%
{\mspace{-2mu}\raisebox{-.6ex}{${\scriptstyle{#1}}$}}}
\nc{\rmats}[1]{{\mathbf{r}}_%
{\mspace{-2mu}\raisebox{-.6ex}{${\scriptscriptstyle{#1}}$}}}
\nc{\shc}{\mathcal{C}}
\nc{\shs}{\mathcal{S}}
\nc{\Fct}{{\on{Fct}}}
\nc{\tC}{\widetilde{\shc}}
\nc{\Zp}{\Z_{\ge0}}
\nc{\tPhi}{\widetilde{\Phi}}
\nc{\tT}{{\widetilde{\T}}}
\nc{\Ob}{\on{Ob}}
\nc{\bwr}{\mbox{\large$\wr$}}
\nc{\Img}{\on{Im}}
\nc{\Ab}{\mathcal{A}^{\mathrm{big}}}
\nc{\Sb}{\mathcal{S}^{\mathrm{big}}}
\nc{\As}{\mathcal{A}}
\nc{\Ss}{\mathcal{S}}
\nc{\ntens}{\widetilde{\otimes}}
\nc{\hR}{\widehat{R}}
\nc{\nconv}{\mathop{\mbox{\large $\odot$}}}
\nc{\snconv}{\mbox{\scriptsize$\odot$}}
\nc{\ts}{\tilde{s}}
\nc{\sho}{\mathcal{O}}
\nc{\bc}{\begin{cases}}
\nc{\ec}{\end{cases}}
\nc{\slnh}{{\widehat{\mathfrak{sl}}_N}}
\nc{\UA}{U_q'(\slnh)}
\nc{\cQ}{\mathcal{Q}}
\nc{\Irr}{\mathcal{I}rr}
\nc{\tQ}{\widetilde{\cQ}}
\nc{\bs}{\mathbf{s}}
\nc{\bL}{\mathbb{L}}
\nc{\tg}{\tilde{g}}

\nc{\conv}{\mathbin{\mbox{\large $\circ$}}}
\nc{\Rm}{R^{\mathrm{ren}}}

\nc{\bQ}{\ol{Q}}
\renewcommand{\Im}{\on{Im}}

\nc{\de}{\on{\textfrak{d}}}

\nc{\xmono}{\ar@{>->}}
\nc{\xepi}{\ar@{->>}}
\nc{\db}[1]{\raisebox{-.5ex}[2ex][1.8ex]{$#1$}}
\nc{\wb}[1]{\mbox{$\rule[-1.1ex]{0ex}{2ex}#1$}}
\nc{\univ}{\mathrm{univ}}
\nc{\rM}{{}^*\mspace{-2mu}M}
\nc{\lM}{M^*}
\nc{\uqm}{\uq\smod}
\nc{\tR}{\widetilde{R}_{\gamma,\beta}}
\nc{\tx}{\tilde{x}}
\nc{\bi}{\mathbf{i}}
\nc{\ttau}{\widetilde{\tau}}

\nc{\tEnd}{\on{\widetilde{E}nd}}
\nc{\tHom}{\on{\widetilde{H}om}}

\nc{\K}{{J}}
\nc{\Kex}{{\K}_{\mathrm{ex}}}
\nc{\Kfr}{{\K}_{\mathrm{f\mspace{.01mu}r}}}
\nc{\coro}{\cor}
\nc{\tB}{\widetilde{B}}
\nc{\seed}{\mathscr{S}}
\nc{\up}{\mathrm{up}}
\nc{\bfa}{\mathbf{a}}
\newcommand{\wB}{\widetilde{B}}

\newcommand{\wl}{\mathsf{P}}   % weight lattice
\newcommand{\comp}{\Delta_+}
\newcommand{\comm}{\Delta_-}

\newcommand{\hconv}{\mathbin{\mbox{$\nabla$}}}

\nc{\tensp}{\otimes_{_+}\mspace{-1mu}}
\nc{\tensm}{\otimes_{_-}\mspace{-1mu}}

\newcommand{\ri}{{\mspace{1mu}\rm r}}

\newcommand{\oi}{\overline{\iota}}
\nc{\bg}{{\oi_\g}}

\nc{\An}{A_q(\n)}
\nc{\tEs}{\widetilde{E}^*}
\def\max{{\mathop{\mathrm{max}}}}
\nc{\pn}{p_\n}
\nc{\dP}{\mathrm{E}^*}
\nc{\Up}{U_q^+(\g)}
\nc{\Ag}{A_q(\g)}
\nc{\QA}{\mathbf{A}}
\nc{\Pd}{\wl^+}
\nc{\Po}{\wl}
\nc{\De}[1]{\Delta(#1)}
\nc{\rt}{\ri}
\nc{\prtl}{\rtl_+}
\nc{\nrtl}{\rtl_-}
\nc{\lt}{\mathrm{l}}
\nc{\wtl}{\wt_\lt}
\nc{\wtr}{\wt_\rt}
\nc{\Cmp}{\comp}
\nc{\Cmm}{\comm}
\nc{\Cm}{\Delta}
\nc{\tUq}[1][{\g}]{\tilde{U}_q(#1)}
\nc{\shf}{\mathcal{F}}
\nc{\Bup}{\mathbf{B}^\up}
\nc{\Blow}{\mathbf{B}^\low}
\nc{\bcom}{\begin{comment}}
\nc{\ecom}{\end{comment}}
\nc{\ms}{\mspace}
\renewcommand{\ge}{\geqslant}
\renewcommand{\le}{\leqslant}

\newlength{\mylength}
\begin{document}

\title{Crystal bases and categorifications}

%%%%%%%%%%%%

\author%[M. Kashiwara]
{Masaki Kashiwara}

\address{Research Institute for Mathematical Sciences, Kyoto University,
Kyoto 606-8502, Japan \& Korea Institute for Advanced Study, Seoul 02455, Korea }
\email{masaki@kurims.kyoto-u.ac.jp}

\thanks{This work was supported by Grant-in-Aid for
Scientific Research (B) 22340005, Japan Society for the Promotion of
Science.}

%\date{\today}

\maketitle

\begin{abstract}
This is a survey of the theory of crystal bases, global bases and cluster algebra structure on the quantum coordinate rings.
\end{abstract}

%\section*{Introduction}

\section{Crystal bases}
The notion of {\em quantum groups} 
(or {\em quantized universal enveloping algebras})
was introduced by Drinfeld \cite{D}
and Jimbo \cite{J} around 1985 in order to explain trigonometric
$R$-matrices in 2-dimensional solvable models
in statistical mechanics.
Since then, the quantum group has been one of 
the important tools to describe new symmetries in representation
theory %Mathematical Physics 
and other fields.
\par
The quantum group
$\Uq$ is an algebra over $\C(q)$, in which $q$ is a parameter of temperature in
the 2-dimensional solvable model,
and $q=0$ corresponds to the absolute temperature zero.
The notion of  crystal bases was motivated
by the belief that something extraordinary should happen
at the absolute temperature zero.
In fact, %as we shall explain in this talk,
it turns out that the representations of $\Uq$ have good bases at
$q=0$, which we call {\em crystal bases} (\cite{K0,K.C}).
Here a basis at $q=0$ of a $\C(q)$-vector space $V$ means  a pair
$(L,B)$ consisting of a free module $L$ of $V$
over the ring $A_0\seteq\set{f\in\C(q)}{\text{$f$ is regular at $q=0$}}$
%(called the crystal lattice)
and a basis $B$ of the $\C$-vector space $L/qL$ together with
an isomorphism $\C(q)\tens_{A_0}L\simeq V$.
%Here, $A_0\seteq\set{f\in\C(q)}{\text{$f$ is regular at $q=0$}}$.

The crystal bases have nice properties such as
uniqueness, stability by tensor products, etc.
Moreover the modified actions of the simple root vectors induce
a combinatorial structure on the crystal basis,
called {\em crystal graph}.
This permits us to reduce many problems in the
representation theory
to problems of the combinatorics.

For example, the combinatorics describing
the representation theory of $\gl_n$
can be well understood with crystal bases.
The irreducible representation over $\gl_n$ is parameterized 
by the highest weight,
which corresponds to a Young diagram in the combinatorial 
language.
Then Young tableaux with the given Young diagram
label the crystal basis of the corresponding irreducible representation.
The Littlewood-Richardson rule,
describing the decomposition of
the tensor product of a pair of representations
into irreducible components,
may be also clearly explained by means of crystal bases.
 
The combinatorial description
of crystal bases by Young tableaux is generalized to
other types of simple Lie algebras in \cite{KN}.
Littelmann gave a completely different combinatorial description
of crystal bases (path model)
in \cite{Li}. Lusztig gave also 
bases at $q=0$ from the PBW bases in the ADE case
(\cite{Lus90}).

See \cite{K.B,K.B1} for surveys of crystal bases.

\section{Global bases}
A crystal basis is a basis at $q=0$. 
However, we can extend this basis to the whole $q$-space 
in a unique way
to obtain a true basis of the representation, which we call
the {\em lower global basis} (\cite{K.C}). 
More precisely, the representation $V$ is a $\C(q)$-vector space and it is
equipped a bar-involution
$c\cl V\isoto V$, i.e., a $\C$-linear involution
of $V$ satisfying $c(qv)=q^{-1}c(v)$ for any $v\in V$.
Moreover $V$ has a certain $\C[q^{\pm1}]$-form, 
i.e., a free $\C[q^{\pm1}]$-submodule
$V_{\C[q^{\pm1}]}$ of $V$ such that $\C(q)\tens_{\C[q^{\pm1}]}V_{\C[q^{\pm1}]}\isoto V$.
Then we can prove that the map $L\cap c(L)\cap V_{\C[q^{\pm1}]}\to L/qL$ is an 
isomorphism and the inverse image $\mathbf{B}\subset L\cap c(L)\cap V_\Q$
of $B$ is a basis of the $\C(q)$-vector space $V$ and also
a basis of the $A_0$-module $L$.
We call $\bB$ the {\em lower global basis} of $V$.

The negative half $\Uqm$ of $\Uq$
has also a crystal basis $B(\Uqm)$, 
and it lifts to the lower global basis 
$\Blow$.

It is first introduced by Lusztig under the name of
canonical basis in the ADE cases 
inspired by the work of Ringel (\cite{Ringel})
describing $U_q^-(\g)$ as the Hall algebras
associated with quivers (see \cite{Lus90,Lus93}).
It is shown that the canonical basis and the global basis coincide
(Grojnowski-Lusztig\,\cite{G-L}).

The dual basis of the lower global basis is called the {\em upper global basis}.
The lower global basis $\Blow$ of $\Uqm$ is a basis
of the integral form
$U^-_{\Z[q^{\pm1}]}(\g)$ of the negative half of $U_q(\g)$.
%the universal enveloping algebra $U(\n)$ of the half $\n$ of $\g$ and
The upper global basis $\Bup$ is the basis of $A_q(\n)$, the dual form
of  the integral form $U^-_{\Z[q^{\pm1}]}(\g)$ (\cite{K.G}).
With a canonical symmetric bilinear form on $\Uqm$, the dual form $A_q(\n)$
may be regarded as a subalgebra of $\Uqm$.
However, at $q=1$, the integral form
$U^-_{\Z[q^{\pm1}]}(\g)$  becomes the universal enveloping algebra
$U(\n)$ of the negative half $\n$ of $\g$,
while $A_q(\n)$ becomes the coordinate ring $\C[\n]$ of $\n$.

\section{Quiver Hecke algebras}
The notion of {\em quiver Hecke algebras} (sometimes called the {\em \KLRs})
is introduced independently by Rouquier (\cite{R08, R11}) and Khovanov-Lauda 
(\cite{KL09, KL11}).

It is a family of $\mathbb Z$-graded algebras which {\em categorifies}
the negative half $\Uqm$ of a quantum group. 
More precisely, 
%if $U_q(\g)$ is a quantum group associated with a
%symmetrizable Cartan datum, then 
there exists a  family of algebras
$\{R(n)\}_{n \in \Z_{\ge 0}}$ (quiver Hecke algebras) such that the (split) Grothendieck group $K(R\gproj)$ of
the direct sum $R\gproj\seteq\soplus\nolimits_{n\in\Z_{\ge0}}R(n)\gproj$ of
the categories of finitely generated projective graded $R(n)$-modules
is isomorphic to the integral form
$U^-_{\A}(\g)$ of the negative half of $U_q(\g)$ as a $\Z[q^{\pm1}]$-algebra. 
Note that the Grothendieck group $K(R\gproj)$ has a natural structure of
$\Z[q^{\pm1}]$-algebra, where the
multiplication is induced by the monoidal category structure of
$R\gproj$ given by the convolution product and
the action of $q$ is induced by the
grading shift functor. 
\Rem 
\bnum\item
The quiver Hecke algebra $\{R(n)\}_{n\in\Z_{\ge0}}$ depends
on a base ring $\cor$, an index set $I$ (the index of simple roots) and
a family of polynomials $Q_{ij}(u,v)\in\cor[u,v]$ ($i\not=j\in I$) such that
$Q_{ij}(u,v)=Q_{ji}(v,u)$.
It is related with $U_q(\g)$ as follows. Let
$A=\{a_{ij}\}_{i,j\in I}$ be a generalized Cartan matrix for $\g$. Then
$Q_{ij}(u,v)$ is a polynomial in $u$ (with coefficients in $\cor[v]$)
of degree $-a_{ij}$ with an element of $\cor^\times$ as the top coefficient. 
(For the precise definition,  see e.g., \cite{R08,KKKO18}.)
\item
For a given generalized Cartan matrix $A$, the representation theory
of $R$ depends on a choice of $\{Q_{ij}\}$ (see \cite{K.P}).
\item
The quiver Hecke algebra is called {\em symmetric} if
$Q_{ij}(u,v)=(u-v)^{-a_{ij}}$ up to a constant multiple.
\ee
\enrem

The Grothendieck group $K(R\gproj)$ has 
the basis consisting of
the isomorphism classes of indecomposable project modules.
If $R$ is symmetric (and the base field is characteristic $0$),
then the lower global basis corresponds to
this basis of $K(R\gproj)$
(Varagnolo-Vasserot \cite{VV09}).

\medskip
There is a dual statement.
The Grothendieck group $K(R\gmod)$ of 
the direct sum $R\gmod\seteq\soplus\nolimits_{n\in\Z_{\ge0}}R(n)\gmod$ of
the categories of finite-dimensional graded $R(n)$-modules
is isomorphic to the dual of $K(R\gproj)$.
Hence $K(R\gmod)$ is isomorphic to $A_q(\n)$, the dual form
of  the integral form
$U^-_{\A}(\g)$.
Moreover, if $R$ is symmetric (and the base field is characteristic $0$),
then the upper global basis corresponds to
the set of the isomorphism classes of simple modules.

As a consequence, the coefficients $C_{b_1,b_2}^{\ b}$ 
appearing in the multiplication 
$$b_1\cdot b_2=\sum\limits_bC_{b_1,b_2}^{\ b}\,b\quad\text{where $b_1,b_2,b\in\Blow$
or $\Bup$}$$
belong to $\Z_{\ge0}[q,q^{-1}]$
when the generalized Cartan matrix is symmetric.
Here, $\Z_{\ge0}[q,q^{-1}]$ is
the set of Laurent polynomials in $q$ with non-negative 
integers as coefficients.
In general, the positivity fails 
(first observed in the G${}_2$ case by S.\ Yamane (\cite{Y}).

Note that Lauda-Vazirani (\cite{LV11}) proved that
the set of the isomorphism classes of simple modules in $R\gmod$
is canonically isomorphic to
$B(\Uqm)$ for an arbitrary quiver Hecke algebra.

The cyclotomic
quotient $R^{\Lambda}(n)$ of $R(n)$ provides a categorification of the
integrable highest weight module $V(\Lambda)$ of $U_q(\g)$
(\cite{KK12}) with a dominant integral weight $\La$ as a highest weight.

\medskip
One of the motivations of these categorification theorems originated
from the so-called {\em LLT-Ariki theory}. In 1996,
Lascoux-Leclerc-Thibon (\cite{LLT}) conjectured that the irreducible
representations of Hecke algebras of type $A$ are controlled by the
upper global basis
 of the basic representation of the quantum affine algebra $U_q(A^{(1)}_{N-1})$.
Soon after, Ariki (\cite{Ariki})  proved this conjecture by showing that the
cyclotomic quotients of affine Hecke algebras categorify the
irreducible highest weight modules over $U(A^{(1)}_{N-1})$, the
universal enveloping algebra of
 affine Kac-Moody algebra of type $A^{(1)}_{N-1}$.
 In \cite{BK09, R08}, Brundan-Kleshchev (\cite{BK09})
and Rouquier (\cite{R11}) showed that
the affine Hecke algebra of type $A$ is isomorphic to 
the quiver Hecke algebra of
type $A^{(1)}_{N-1}$ or of type $A_\infty$ up to a specialization and a
localization. Thus the quiver Hecke algebras
can be understood as a graded version
of the affine Hecke algebras of type $A$, and
Kang-Kashiwara's
cyclotomic categorification theorem (\cite{KK12}) is a generalization of Ariki's
theorem on type $A^{(1)}_{N-1}$ and $A_\infty$ to all symmetrizable Cartan
datum.

\section{Cluster algebras}
As one can imagine from the fact that
$A_q(\n)$ is a commutative algebra at $q=1$,
the upper global basis of $A_q(\n)$ has an interesting multiplicative property.

Berenstein and Zelevinsky (cf.\ \cite{BZ})
conjectured that, when the the generalized Cartan matrix is of finite type,
there exists a family $\shf$ of
finite subsets of the upper global basis $\Bup$ of
$A_q(\n)$ satisfying the following properties:
\bna
\item 
Any pair $(x,y)$ of elements of $C\in\shf$ is $q$-commutative
(i.e., $xy=q^nyx$ for some $n\in\Z$),
\item For any $C\in\shf$, 
any $C$-monomial, i.e., an element of the form
$x_1\cdots x_\ell$ with $x_1,\ldots,x_\ell\in C$, 
belongs to $q^\Z\Bup\seteq\set{q^nb}{n\in\Z,b\in\Bup}$.
\item $\Bup$ is the union of $C$-monomials (up to a constant multiple)
where $C$ ranges over $\shf$.
\ee
Then Leclerc (\cite{L03}) gave a counterexample to this conjecture.
He called an element $b\in\Bup$ a {\em real} vector if
$b^2\in q^\Z\Bup$. Otherwise, $b$ is called {\em imaginary}.
If Conjecture were true, then any $b\in\Bup$  would be a real vector.
He gave examples of imaginary vectors
%for all finite-dimensional simple Lie algebras
%except type $A_n$ ($n\le 4$) and $B_2=C_2$.
for types A${}_n$ ($n\ge5$), B${}_n$ ($n\ge3$), C${}_n$ ($n\ge3$),
D${}_n$ ($n\ge4$), and all exceptional types.

Although their conjecture failed, 
their idea survives and it was one of the motivations
of the introduction of cluster algebras
by Fomin and Zelevinsky (\cite{FZ02}).
They replaced condition (c) with 
a weaker condition

\be
\item[{(c)${}^\prime$}]
$\Uqm$ is generated by $\bigcup_{C\in\shf}C$ as a $\C(q)$-algebra.
\ee
They call $C\in\shf$ a cluster,
and reformulated the conjecture in the language of cluster algebras.

%and the clusters are related to each other
%by `` mutations''.

\medskip

A {\em cluster algebra} is a $\Z$-subalgebra of a rational function field
given by a set of generators,
called the {\it cluster variables}. These generators are grouped
into overlapping subsets, called {\it clusters}.
The clusters are defined
inductively by a procedure called {\it mutation} from the {\it
initial cluster} $\{ X_i\}_{1 \le i \le r}$.
The mutation is controlled by
an exchange matrix $\wB=(b_{ij})_{ij}$ as follows. 
By the mutation at $k$ ($1\le k\le r$),
a new cluster is created from the old cluster by replacing the $k$-th variable $X_k$ with
\eq
X'_k=\dfrac{\prod\limits_{i;\;b_{ik}>0}X_i^{b_{ik}}+\prod\limits_{i;\;b_{ik}<0}X_i^{-b_{ik}}}
{X_k}\;.\label{eq:mutation}
\eneq
The exchange matrix $\wB=(b_{ij})_{ij}$ is also mutated to
the new exchange matrix $\mu_k(\wB)=(b'_{ij})_{ij}$
given by
\eqn
b'_{ij}&=&\bc
-b_{ij}&\text{if $i=k$ or $j=k$,}\\
b_{ij}+(-1)^{\delta(b_{ik}<0)}\max(b_{ik}b_{kj},0)\hs{3ex}&\text{otherwise.}
\ec
\eneqn
A {\it cluster monomial} is a monomial of cluster
variables in one cluster. 
In the conjecture above, every $C\in\shf$ is a cluster.
Moreover, any two members of $\shf$ 
are connected by successive mutations.
We do not assume (c) but we assume that the algebra is generated 
(as an algebra) by the cluster monomials.

Fomin and Zelevinsky proved that every cluster variable is a Laurent
polynomial of the initial cluster $\{ X_i\}_{1 \le i \le r}$.
They conjectured that this Laurent polynomial has positive
coefficients (\cite{FZ02}). This {\it positivity conjecture} was
proved by Lee and Schiffler in the {\it skew-symmetric} cluster algebra case
in \cite{LS13}.
The {\it linearly independence
conjecture} on cluster monomials was proved
in the skew-symmetric cluster algebra case
in \cite{CKLP12}.

The notion of quantum cluster algebras, introduced by Berenstein and Zelevinsky in
\cite{BZ05}, can be considered as a $q$-analogue of cluster algebras.
It is an algebra over $\Z[q^{\pm1/2}]$.
The cluster variables in a cluster $q$-commute with each other.
%a skew-symmetric matrix $L$.
As in the cluster algebra case, every cluster variable belongs to
$\Z[q^{\pm 1/2}][X_i^{\pm 1}]_{1 \le i \le r}$ 
for the initial cluster $\{ X_i\}_{1 \le i \le r}$
(\cite{BZ05}), and is
expected to be an element of $\Z_{\ge0}[q^{\pm 1/2}][X_i^{\pm 1}]_{1 \le i \le r}$,
which is referred to as the {\it quantum
positivity conjecture} (cf.\ \cite[Conjecture 4.7]{DMSS}).
In \cite{KQ14}, Kimura and Qin proved the quantum positivity conjecture for quantum cluster algebras containing {\it acyclic} seed
and specific coefficients.

Assume that the generalized Cartan matrix is symmetric.
In a series of papers \cite{GLS11,GLS05,GLS}, Gei\ss, Leclerc and
Schr{\"o}er showed that the quantum unipotent coordinate algebra
$A_q(\mathfrak{n}(w))$
has a
skew-symmetric quantum cluster algebra structure whose initial cluster
consists of {\it quantum minors}.
Here, $A_q(\mathfrak{n}(w))$ is a $\Z[q^{\pm1}]$-subalgebra
of $A_q(\mathfrak{n})$ associated with a Weyl group element $w$.
In \cite{Kimu12}, Kimura
proved that $A_q(\mathfrak{n}(w))$ is {\it compatible} with the
upper global basis $\B^{\upper}$ of $A_q(\mathfrak{n})$; i.e., the
set $\B^{\upper}(w) \seteq A_q(\mathfrak{n}(w)) \cap \B^{\upper}$
is a basis of $A_q(\mathfrak{n}(w))$. Thus, %with a result of \cite{CKLP12}, 
one can expect that every cluster monomial
of $A_q(\mathfrak{n}(w))$ is contained in the upper global basis
$\B^{\upper}(w)$, which is named {\it the quantization conjecture} by Kimura
(\cite{Kimu12}).
\section{Monoidal categorification}
This conjecture (in the symmetric generalized Cartan matrix case)
is proved affirmatively
by Kang-Kashiwara-Kim-Oh(\cite{KKKO18}), using the
monoidal categorification of $A_q(\mathfrak{n}(w))$
by a subcategory of 
the module category $R\gmod$ over the quiver Hecke algebras.

%In particular any cluster monomial is a real element of the upper global basis.

%\begin{comment}
In \cite{HL10}, Hernandez and Leclerc introduced the notion of {\em
a monoidal categorification of  a cluster algebra}.
Let $(\shc,\,\tens)$ be a monoidal category.
We say that a simple object $S$ of $\shc$ 
is {\em real} if $S \tens S$ is simple. 
We say that a simple object $S$ is {\em prime} if there exists no non-trivial
factorization $S \simeq S_1 \tens S_2$. 
They say that $\shc$ is
a monoidal categorification of a cluster algebra $A$ if
the Grothendieck ring of $\shc$ is isomorphic to
$A$ and if

\vs{1.5ex}
\hs{0ex}\parbox{80ex}{

\begin{enumerate}
\setlength{\itemsep}{3pt}
\item[{\rm (M1)}] any cluster monomial of $A$ corresponds to
the class of a real simple object of $\shc$,
\item[{\rm (M2)}] any cluster variable of $A$ corresponds to
the class of 
a real simple prime object of $\shc$.
\end{enumerate}}

\vs{1.5ex}
\noi
(Note that the above version is
weaker than the original definition of the monoidal categorification in
\cite{HL10}.) 
They proved that certain categories
of modules over symmetric quantum affine algebras
$U_q'(\g)$ give monoidal categorifications of cluster algebras.
Nakajima extended it to the cases of the cluster algebras of type $A,D,E$
(\cite{Nak11}) (see also \cite{HL13}).

Once a cluster algebra $A$ has a monoidal categorification,
the positivity of cluster
variables of $A$ and the linear 
independence of cluster monomials of $A$ follow.

%(see \cite[Proposition 2.2]{HL10}).

%In this paper, we will refine their notion of
%monoidal categorifications including the quantum cluster algebra case.

%
%\medskip
%
%The Khovanov-Lauda-Rouquier (or simply KLR) algebras, introduced by
%Khovanov-Lauda \cite{KL09,KL11} and Rouquier \cite{R08}
%independently, are a family of $\Z$-graded algebras which
%categorifies the negative half $U_q^-(\g)$ of a {\it symmetrizable}
%quantum group $U_q(\g)$. More precisely, there exists a family of
%algebras $\{ R(-\beta) \}_{\beta \in \rtl^-}$ such that the Grothendieck
%ring of $R \gmod \seteq \bigoplus_{\beta \in \rtl^-}R(-\beta)\gmod$, the direct sum
%of the categories of finite-dimensional graded $R(-\beta)$-modules, is
%isomorphic to the integral form $A_q(\mathfrak{n})_{\Z[q^{\pm1}]}$ of
%$A_q(\mathfrak{n}) \simeq U_q^-(\g)$. Here the tensor functor $\tens$
%of the monoidal category $R \gmod$ is given by
%the convolution product $\conv$, and the action of $q$ is given by
%the grading shift functor. In \cite{VV09, R11},
%Varagnolo-Vasserot and Rouquier
%proved that the upper global basis $\B^\upper$ of $A_q(\mathfrak{n})$
%corresponds to the
%set of the classes of all {\it self-dual} simple modules of $R
%\gmod$ under the assumption that $R$ is associated with a {\it
%symmetric} quantum group $U_q(\g)$.

In order to give a monoidal categorification of
$A_q(\mathfrak{n}(w))$, we use a monoidal cluster,
a categorification of a cluster.
It is a finite set of real simple objects $\{M_i\}_{1\le i\le r}$
in the monoidal category $(R\gmod,\,\tens)$ of the 
finite-dimensional graded modules
over the quiver Hecke algebras,
which satisfies the condition: 
$M_i\tens M_j\simeq M_j\tens M_i$ up to a grading shift.
Then the mutation at $k$ ($1\le k\le r$)
creates a new monoidal cluster 
from the old monoidal cluster by replacing the $k$-th object $M_k$ with
a real simple object $M'_k$.
Here, the mutated object $M'_k$ is related with
the original monoidal cluster
by the following exact sequences
(up to grading shifts) instead of the relation \eqref{eq:mutation}:
\eq
&&\ba{l}
0\To\bigotimes\limits_{i;\;b_{ik}<0}M_i^{\otimes\,(-b_{ik})}\To M_k\tens M'_k
\To\bigotimes\limits_{i;\;b_{ik}>0}M_i^{\otimes\, b_{ik}}\To0,\\[2ex]
0\To\bigotimes\limits_{i;\;b_{ik}>0}M_i^{\otimes\, b_{ik}}\To M'_k\tens M_k
\To\bigotimes\limits_{i;\;b_{ik}<0}M_i^{\otimes\,(-b_{ik})}\To0.
\ea
\label{eq:monoidalmutation}
\eneq
Note that if one passes to the Grothendieck group level,
then \eqref{eq:mutation} is a consequence of \eqref{eq:monoidalmutation}.
% implies 

\smallskip
In \cite{KKKO18}, it is proved that if the first %and the second 
step mutations
starting from the initial monoidal cluster
are possible,
then every successive mutations are possible.
As its consequence, %we prove that 
the quantum unipotent coordinate algebra
$A_q(\mathfrak{n}(w))$ has a monoidal categorification, and
the quantization conjecture follows.
Note that 
F.\ Qin also provided a proof of the conjecture for a large class with
a condition on the Weyl group element $w$  in a completely different method
(\cite{Qin15}).

Note that the converse of (M1)
is conjectured by several experts (e.g., see \cite{GLS11,Kimu12}):

\begin{enumerate}
\item[{\rm (M1)${}^\prime$}] 
%the classes of 
any real element of $B\bl A_q(\mathfrak{n}(w))\br$
is a cluster monomial.
\end{enumerate}

It is still open.

%is proved by the monoidal categorification.
%As a byproduct,
%it is also proved that 
%any cluster monomial is a real element of the upper global basis.

%Combining works of \cite{GLS,Kimu12,VV09}, the quantum unipotent
%coordinate algebra $A_q(\mathfrak{n}(w))$ associated with a
%symmetric quantum group $U_q(\g)$ and a Weyl group element $w$
%is isomorphic to the Grothendieck group of
%the monoidal abelian full subcategory $\shc_w$ of $R \gmod$
%satisfying the following properties: {\rm (i)} $\shc_w$ is
%stable under extensions and grading shift functor, {\rm (ii)} the
%composition factors of $M \in \shc_w$ are contained
%in $\B^{\upper}(w)$. However it is not evident the conditions (M1) and (M2)
%are satisfied. The purpose of this paper is provide a theoretical background
%in order to prove (M1) and (M2). In the forthcoming paper, we
%prove (M1) and (M2) as an application of the results of the present paper.

\section{Real objects}
%A simple object $M$ of an abelian monoidal category $\shc$ is called real
%if $M\tens M$ is also simple.
By the monoidal categorification,
a real element of the upper global basis corresponds to a real object
of the monoidal category $R\gmod$ of finite-dimensional graded modules
over the quiver Hecke algebra.
Real simple objects have remarkable properties.
Let $(\shc,\tens)$ be a monoidal category which is either
the monoidal category $R\gmod$ of finite-dimensional graded modules
over the {\em symmetric} quiver Hecke algebra,
% of finite-dimensional graded modules over the quiver Hecke algebra
or the monoidal category of finite-dimensional modules
 over an affine quantum group.
Let $\cor$ be the base field of $\shc$.
Then we have the following propositions.
\Prop[{\cite{KKKO14, KKKO18}}]
Let $M$ and $N$ be simple objects of $\shc$.
We assume that either $M$ or $N$ is real.
Then we have \ro forgetting grading shifts\rf
\bnum
\item $\Hom(M\tens N,M\tens N)=\cor\id_{M\tens N}$,
\item
there exists a non-zero morphism $\rmat{}\cl M\tens N\to N\tens M$ such that
$$\Hom(M\tens N,N\tens M)=\cor\ms{1.5mu}\rmat{},$$
\item
$\Im(\rmat{})$ is simple and it coincides with the head of $M\tens N$
and also with the socle of $N\tens M$.
\ee
\enprop
Conversely, a simple object $M$ is real as soon as
$\End(M\tens M)=\cor\id_{M\tens M}$.

We denote by $M\hconv N$ the simple head of
$M\tens N$.

\Prop[\cite{KKKO14}]
Let $M$ be a real simple object of $\shc$ and let $N,N'$ be simple objects.
%\bnum
%\item
If $M\hconv N\simeq M\hconv N'$, then $N\simeq N'$.
%\item If there is a non-zero morphism
%$M\tens N\to N'\tens M$,
%\ee
\enprop

\Rem
\bnum
\item
Note that the Grothendieck ring  $K(\shc)$
is commutative if $\shc$ is the module category over an affine quantum group.
Similarly, the Grothendieck ring  $K(R\gmod)$
is commutative if we forget the grading shifts.
\item
In \eqref{eq:monoidalmutation}, we have
$$\text{$M_k\hconv M'_k\simeq
\bigotimes\limits_{i;\;b_{ik}>0}M_i^{\otimes\, b_{ik}}$ and
$M'_k\hconv M_k\simeq
\bigotimes\limits_{i;\;b_{ik}<0}M_i^{\otimes\, -b_{ik}}$.}$$
\ee
\enrem

\end{document}